\documentclass{article}

\usepackage{times}
\usepackage{latexsym}
\usepackage{amsmath}   
\usepackage{amsfonts}  
\usepackage{amssymb}   
\usepackage{amsthm}    
\usepackage{epsfig}    
\usepackage{psfrag}
\usepackage{bm}
\usepackage{graphics}
\usepackage{graphicx}
\usepackage{pstricks}
\usepackage{pst-node}

\newcommand{\pr}{\mbox{\sf P}}
\newcommand{\ex}{{\bf\sf E}}               

\newcommand{\al}{\alpha}                
\newcommand{\G}{\Lambda}                

\newtheorem{thm}{Theorem}
\newtheorem{lem}[thm]{Lemma}

\newtheorem{defn}[thm]{Definition}

\def\th{\theta}

\begin{document}

\title{Stationary Performance  Analysis of Grishechkin
  Processor-Sharing Queues: An Integral
  Equation Approach
}


\author{Yingdong Lu \\ IBM T.J. Watson Research Center \\ Yorktown Heights, NY 10598}



\date{March 2021}

\maketitle

\begin{abstract}

\vskip 0.2cm

We compute the stationary performance metrics of a single server $M^X/G/1$ queue under a class of generalized processor-sharing scheduling
policies that are proposed by Grishechkin. This class of processor-sharing policies allow service
capacities to be allocated to jobs based on the amount of service they
attained. In \cite{grishechkin}, Grishechkin derives an integral equation that is satisfied by the Laplace transform of the stationary performance metrics under these policies. Our main focus in this paper is to derive the solution to this integral equation. This is achieved through a series of transforms that convert the integral equation into a more tractable form, then solve it. Then we derive approximations to the density functions of the performance metrics through inverting the Laplace transforms. In the last part of the paper, we apply our results to some well-known scheduling policies that are either a special case of the Grishechkin processor-sharing policies or can be treated as the limits of a sequence Grishechkin processor-sharing policies. These examples include the egalitarian processor-sharing policy, the discriminatory processor-sharing policy, shortest residual processing time first rule, foreground-background policy and the time function scheduling policy.

\end{abstract}

\noindent
{\bf Keywords: } {Queueing system; Laplace transform; Grishechkin processor-sharing}

\section{Introduction}
\label{sec:intro}

\vskip 0.2cm

In this paper, we consider a $M^X/G/1$ queue, that is, a single server
queue with Poisson batch arrivals and general service
distribution. We assume that the service distribution could follow any
probability distribution with finite support. Therefore, the queueing
model considered in the paper is fairly general and flexible. Our goal is to
obtain an analytic characterization of its stationary performance
metrics, such as the stationary system size under a class of general processor-sharing scheduling
policies.

There are many reasons for studying a queueing model with
a processor-sharing scheduling policy, which has been
proved to be a powerful and flexible model in performance
analysis. First, it has been demonstrated through extensive
research from different aspects that a processor-sharing queue is an accurate mathematical model for the performance of a real life system that operates under scheduling policies such as round-robin or time-sharing,
and these policies are commonly adapted in job scheduling for computer
or communications networks, see, e.g. \cite{Kleinrock}. Secondly, in
some recent study of fairness of the scheduling policies,
processor-sharing policy often serves as a reference or a base on which other
more complicated policies are developed and analyzed, see,
e.g. \cite{Wierman,FriedmanHenderson}. Fairness is a critical issue
for many new applications,  especially those in web services. These
facts all contribute to the importance of the performance analysis of
processor-sharing policies.

One important feature of processor-sharing is its flexibility. The basic
principle of processor-sharing is to allocate resource capacity
proportional to jobs that need to be processed. However, there can be
many different ways to decide the proportion for each job, besides the equal proportion that corresponding to the
egalitarian processor-sharing policy. Different regimes for
determining the proportion can produce very different performances,
which have to be considered very carefully when a scheduler is designed
for any types of systems. A very natural scheme is to let the
proportions be determined by the service already attained by the
jobs, because this is usually a quantity that can be obtained with low
overhead charged to many systems.

The success of the deployment of these policies and ideas in the
design and performance analysis of the systems depends on efficient
evaluation of a wide class of processor-sharing polices. Meanwhile,
the theoretical development for understanding a processor-sharing
queue still has plenty of room for
improvements. It starts in the case of systems
with memory-less assumptions on distributions (Poisson arrival or
exponential service time), see, for example, Yashkov \cite{yashkov87}
for a survey on results of this nature, as well as other developments
prior to $1987$. For a general $GI/GI/1$ queue, Jean-Marie and Robert
\cite{jean-marie} study the transient behavior  to obtain the rate of
growth of the number of jobs in the system; Baccelli and Towsley
\cite{baccelli} investigate stochastic ordering properties of the
response time; Grishechkin \cite{grishechkin91a, grishechkin}
establishes a strong approximation limit for the number of jobs in the
system under heavy traffic; More recently, Gromoll {\it et al.}
\cite{gromoll2002, gromoll2004} establish fluid and diffusion limits
through measure-valued processes.

In his effort to obtain the heavy traffic limit of the steady state
distribution of the performance metrics, Grishechkin establishes a
connection between the state process of a $M^X/G/1$ queue under a very general class of processor-sharing policies and a Crump-Mode-Jagers branching process, then derives an
integral equation for a wide class of performance metrics from the
dynamics of the branching process. Grishechkin names this class of policies generalized processor-sharing policies. To avoid any
confusion with the generalized processor sharing policy for multiclass
networks, we will call this class of policies Grishechkin
processor-sharing policy hereafter. Under a  Grishechkin
processor-sharing policy, each job in the
system receives a portion of the processing capacity, and the
proportion is determined by a random function whose variables are the
amount of service each job has already received up to the moment under
consideration. By choosing different random functions,  Grishechkin
processor-sharing can cover a
wide range of known scheduling policies. For example, the egalitarian processor-sharing policy is a special case of Grishechkin
processor-sharing when the random function just assigns the same weight to each job in the system, and the policy of shortest residual processing time first rule can be treated as the limit of a sequence of  Grishechkin processor-sharing policies with different parameters.

Due to the complicated nature of this integral equation, only the
simplest case of egalitarian processor-sharing with unit demand has been solved in
\cite{grishechkin}. In this paper, our main contribution is to discover that several transforms can be applied to the integral equation so that it can be reduced to a
more tractable form. To be more precise, we reduce the integral
equation to a special form of combined first and second kind, and then
we observe some special structures of the classical solutions to the first
and second kinds of integral equation with the kernel appeared in our integral equation. These observations enable us to solve the
transformed integral equation to obtain the Laplace transform of the
stationary performance metrics. In the current paper, we will focus on
one performance metric, the system size, that is, the number of jobs
in the system. Similar approaches can be applied to other metrics that
are also studied in \cite{grishechkin91a, grishechkin}.

Then we devote our efforts to inverting the Laplace transform we obtained through the integral equation. Because in the course of solving the integral equation, we need to make use of the Laplace transform solution to the classical integral equation of the first kind, the final Laplace transform is essentially a two-dimensional Laplace transform. We adapt sophisticated Laplace transform inverting method to our problem, and obtain approximations to the stationary distribution of the performance metrics with accuracy guarantees.

The rest of the paper will be organized as follows. In
Sec. \ref{sec:prelim}, we will give a detailed description of the  Grishechkin
processor-sharing policies, as well as the integral equation for the stationary
system size. In Sec. \ref{sec:CMJ}, we list some basic facts on
the Crump-Mode-Jagers branching process and describe Grishechkin's
results for connecting the queueing process and the branching process.
In Sec. \ref{sec:inteqn}, we will derive the solution for
the integral equation. In Sec. \ref{sec:performance} we will summarize
our results with a final expression of the stationary system size, and
derive numerical procedures for several special cases. In
Sec. \ref{sec:conclusions}, we will conclude the paper with a summary of
our results and contributions, as well as a description of some
ongoing research.

\section{Models and Preliminaries}
\label{sec:prelim}

\vskip 0.2cm

Let us describe the queueing system in detail, and the Grishechkin
processor-sharing scheduling
policies in particular. Jobs arrive at the queueing system following a
compound Poisson process with rate $\Lambda$ and batch size
distribution $B$ with finite support. For each job, the amount of service required is independent and
follows an identical general distribution, let us denote it as $\ell$, and there exists $M>0$ such that $\pr[\ell>M]=0$. We assume that the
arrival and service processes are independent, the inter-arrival and
the batch distributions for different jobs are also independent. The
server serves at the rate of one, and can serve multiple jobs
simultaneously. At time $T=0$, we assume that there are $N$ jobs in
the system, and they attain no service before time $T=0$.

Each job $n$, $n=1, 2,\ldots, $ is associated with a pair of
independent stochastic processes, $(A_n(\cdot), D_n(\cdot))$, indexed
by the amount of service it received. These pairs are i.i.d. for all
jobs and a generic member $(A(\cdot),D(\cdot))$ satisfies,
\begin{itemize}
\item
$A(\cdot)$ and $D(\cdot)$ have absolutely continuous sample path almost surely;
\item
The random variable $\ell = \min\{ T>0, A(T) =0\}$ is almost surely
finite;
\item
The integral $\int_0^\ell (A(y))^{-1}dy $ converges almost surely.
\end{itemize}

Intuitively, the process $A(\cdot)$ reflects the weight of the job,
and $D(\cdot)$ the quantity of interest. Define $Q(t)$ as the number
of jobs in the system at any time $t$, and $V_i(t)$ as the amount of
service attained by each job $i, i=1,2,\ldots, Q(t)$, up to time $t$.
Then policy will allocate service capacity among jobs in the system
such that the following relationship for the service rates of job
$i$ must satisfy,
\begin{eqnarray}
\label{eqn:service_rate}
\frac{d V_i(t) }{dt } = \frac{A_i(V_i)}{\sum_{j=1}^{Q(t)}A_j(V_j)}.
\end{eqnarray}
That is, the sharing of the service capacity is determined by the
distribution of $A$ and the amount of service attained. Define
\begin{eqnarray*}
Q^D(T) =\sum_{i, T_i\le T} D_i (V_i(T)),
\end{eqnarray*}
where $T_i$ denotes the arrival time of job $i$. Then the main
performance metric of interest is the stationary distribution of the
above random variable. Depending upon the functional form of $D(\cdot)$,
$Q^D(T)$ corresponds to different types of performance metrics at
time $T$,
\begin{itemize}
\item
$D(y) = {\bf 1}\{y \in[0, \ell]\}$, $Q^D(T)$  represents the system size;
\item
$D(y) = {\bf 1}\{ y \in[0, \min[a, \ell]]\}$, $Q^D(T)$  represents
  number of jobs with attained service less than $a$;
\end{itemize}
where ${\bf 1}A$ denotes the indicator function of set $A$, and $Q^D$
denotes its stationary distribution. While the methods developed in this paper certainly are not restricted to, for the ease of exposition, we will only present the case when $D(y) = {\bf 1}\{y \in[0, \ell]\}$, i.e. the case in which $Q^D(T)$ represents the system size, or equivalently, the number of jobs that are in the system, so, unless otherwise noted, we will omit the superscript $D$.

\section{Crump-Mode-Jagers branching process and the derivation of the main integral equation}
\label{sec:CMJ}

\vskip 0.2cm

It is shown in \cite{grishechkin}, see Theorem 2.1 in \cite{grishechkin}, that system size process will have the same
probability law as that of the population of a Crump-Mode-Jagers
branching process. In this section, for the completeness of the arguments, we will
introduce some basic facts about the Crump-Mode-Jagers branching
process and sketch the main arguments in \cite{grishechkin91a} for the
derivation of the integral equations that the stationary performance
metric will satisfy.

The following theorem of Grishechkin provides a description of the 
stationary distribution of the system size.

\begin{thm}
The function $\phi(t;u, v)$ satisfies the following integral equation,
\begin{eqnarray}
\label{eqn1_branch_2} &&\phi(t;u, v) = \ex\left[\exp\left\{h(t;u,v)+ \G
\int_0^t f(\phi(t-y;v,u)-1) C(y)dy \right\}\right].
\end{eqnarray}
where $$h(t;u,v) := u \chi(t) + vC(t).$$
\end{thm}
Our performance analysis is reduced to solving this integral
equation.

\section{Solving the Integral Equation}
\label{sec:inteqn}

\vskip 0.2cm

It is well-known that although very common in various fields of
applied mathematics and engineering, an integral equation is rather
difficult to solve, unless it is in the standard first or second kind
integral equation forms, see, e.g. \cite{PorterStirling}. Our integral
equation (\ref{eqn1_branch_2}) does not appear to fall into those
categories. In the following, we will take two steps in obtaining its
solution. First, a sequence of transformations are applied to
(\ref{eqn1_branch_2}) aimed at simplifying it to a more tractable
form. As a result, we found that (\ref{eqn1_branch_2}) can be reduced
to the form of a combined first and second kinds of integral equation.
Then, in the second step, we derive the solution of the combined
first and second kinds of integral equation solution making use of the
special function form of the solution. There is no systematic approach
to solve a combined first and second kinds of integral equation
solution, to the best of our knowledge, our solution is among very few
cases where an explicit solution has been derived.

First, the assumption of the absolute continuity of the $(A(\cdot),
D(\cdot))$ enables us to write the integration equation as,
\begin{eqnarray*}
\phi(t;u,v) &= &\int_{{\mathbb R}^3\times [0,t]}
  \exp[h_1(t;u,v;z_1,z_2,z_3, z_4)]
  \times \nonumber \\ && \exp \left\{
\int_0^t [f(\phi(t-y;u,v))-1] h_2(y, z_1,z_2,z_3, z_4)dy\right\}\\
  &&g(z_1,z_2,z_3,z_4)dz_1dz_2dz_3dz_4.
\end{eqnarray*}
where $g(z_1,z_2,z_3,z_4)$ is the joint density function of $(A(\cdot),
D(\cdot), \ell)$ at time $z_4$, and
\begin{eqnarray}
\label{eqn:H1}
h_1(t;u,v;z_1,z_2,z_3,z_4)&=&[u\chi(z_1,z_2,z_3)(t) +
  vC(z_1,z_2,z_3)(t)]{\bf 1}\{t=z_4\},
\end{eqnarray}
and
\begin{eqnarray}
\label{eqn:H2}
h_2(t,z_1,z_2,z_3,z_4)=C(z_1,z_2,z_3)(t){\bf 1}\{t=z_4\}.
\end{eqnarray}
In the rest of the section, when transforms are applied to the
integral equations, all the actions on $z_1, z_2, z_3$ and $z_4$, in
fact, will be identical, meanwhile our results are not restricted by
any specific function form $g(\cdot)$ takes. Therefore, for the ease
of exposition, we will just use one variable $z$ and a univariate
function $g(z)$, with finite support $[0,M]$, to represent them, and restore their original form at
the end of the derivation. Therefore, the integral equation under
consideration will now bear the following form,
\begin{eqnarray}
\label{eqn:inteqn_1}
&&\phi(t;u,v) =  \int_{{\mathbb R}} \exp[h_1(t;u,v;z)] \times
\nonumber \\ && \exp \left\{
\int_0^t [f(\phi(t-y;u,v))-1] h_2(y,z)dy\right\}g(z)dz.
\end{eqnarray}
with $h_1(t;u,v;z)$ and $h_2(t,z)$ in place for functions defined in
(\ref{eqn:H1}) and (\ref{eqn:H2}).

Define an operator $T:C^1({\mathbb R}_+^4) \rightarrow C^1({\mathbb
  R}_+^3) $, where $C^1(\Omega)$ denotes the family of 
functions on domain $\Omega$ with continuous derivatives,  as follows,
\begin{eqnarray*}
T (\psi(t;u,v;z)) = \int_{\mathbb R} \exp[h_1(t;u,v;z)] \psi(t;u,v;z)
g(z) dz.
\end{eqnarray*}
Let us denote, $\psi(t;u,v;w)$ as the solution of the following
integral equation,
\begin{eqnarray}
\label{eqn:inteqn_2}
&& \psi(t;u,v;w) =  \exp\left[\int_0^t -h_2(y,w)dy\right]\nonumber \\ && \exp \left\{
\int_0^t h_2(y,w) f(T (\psi(t-y;u,v;z))dy \right\}.
\end{eqnarray}
Then, it can be verified through direct calculation that,
\begin{lem}
\label{lem:transform}
\begin{eqnarray*}
\phi(t;u,v)=T (\psi(t,u,v,z))
\end{eqnarray*}
is a solution to the original integral equation (\ref{eqn:inteqn_1}).
\end{lem}
\vskip 1cm

Taking logarithms, the integral equation (\ref{eqn:inteqn_2}) can be
written in the following equivalent form,
\begin{eqnarray*}
&& \log(\psi(t;u,v;w))+ \int_0^t h_2(y,w)dy =\int_0^t
  h_2(y,w) f(T (\psi(t-y;u,v;z))dy.
\end{eqnarray*}
Thus, we can further reduce the problem through the following lemma.
\begin{lem}
\label{lem:log_trans}
Suppose that $\Psi(t,u,v,w_1, w_2, \ldots, w_n)$ is the solution to
\begin{eqnarray}
\label{eqn:inteqn_n}
&&\log(\Psi(t,u,v,w_1, w_2, \ldots, w_n))+  \int_0^t
\sum_{i=1}^nh_2(y,w_i)dy \\ &=&\int_0^t
\sum_{i=1}^nh_2(y,w_i)\left[\int_{{\mathbb R}^n}
\sum_{\ell=0}^n f_\ell \exp\left(\sum_{k=1}^\ell h_1(t;u,v,z_k)\right)
\right. \nonumber \\ && \left.\prod_{k=\ell+1}^n \nu(z_k)\prod_{k=1}^n
g(z_k)\Psi(t,u,v,z_1, z_2, \ldots, z_n)\right. dz_1dz_2\ldots dz_n\Big]dy ,
\nonumber
\end{eqnarray}
where $\nu$ is any probability measure that has support on $[M,
  \infty)$, then,
\begin{eqnarray*}
\phi(t;v,u)=T \left((\Psi(t,u,v, z, z, \ldots, z)^{\frac{1}{n}}\right)
\end{eqnarray*}
is the solution to the original integral equation (\ref{eqn:inteqn_1}).
\end{lem}
{\bf Proof }
First, let us consider the special case of $f(x)= x^n$ for some $n$,
we can rewrite the equation as,
\begin{eqnarray*}
&&\log(\psi(t,u,v,w))+ \int_0^t h_2(y,w)dy \\ &=&\int_0^t h_2(y,w)\left[\int_{{\mathbb R}^n}
\exp\left(\sum_{k=1}^nh_1(t;u,v;z_k)\right) \right. \\&& \left. \times
\prod_{k=1}^n g(z_k)\prod_{i=1}^n\psi(t-y,u,v,z_i)dz_1dz_2\ldots
dz_n\right]dy.
\end{eqnarray*}
For each $i=1,2\ldots, n$, define,
\begin{eqnarray*}
&&\log(\psi(t,u,v,w_i)) + \int_0^t h_2(y,w_i)dy    \\ &=&\int_0^t h_2(y,w_i)\left[\int_{{\mathbb R}^n}
\exp\left(\sum_{k=1}^nh_1(t;u,v;z_k)\right)  \right. \\&&
\left. \times \prod_{k=1}^n
g(z_k)\prod_{i=1}^n\psi(t-y,u,v,z_i)dz_1dz_2\ldots dz_n\right]dy.
\end{eqnarray*}
Next, define
\begin{eqnarray*}
\Psi(t,u,v, w_1, w_2, \ldots, w_n)= \prod_{i=1}^n &\psi(t,u,v,w_i).
\end{eqnarray*}
Therefore, $\Psi(t,u,v,w_1, w_2, \ldots, w_n)$ satisfies,
\begin{eqnarray}
\label{eqn:expan11}
&&\log(\Psi(t,u,v,w_1, w_2, \ldots, w_n)) + \int_0^t \sum_{i=1}^nh_2(y,w_i)dy\\
&= & \nonumber
\int_0^t \sum_{i=1}^nh_2(y,w_i) \int_{{\mathbb R}^n}
\exp\left[\sum_{k=1}^nh_1(t;u,v;z_k)\right] \\ \nonumber &&\times\prod_{k=1}^n
g(z_k)\Psi(t-y,u,v,z_1,z_2, \ldots, z_n)dz_1dz_2\ldots dz_ndy.
\end{eqnarray}
Define,
\begin{eqnarray*}
h_4(t;u,v;  z_1,z_2, \ldots, z_n) =
\exp\left[\sum_{k=1}^nh_1(t;u,v;z_k)\right] \prod_{k=1}^n g(z_k).
\end{eqnarray*}
Then the above expression can be written as,
\begin{eqnarray*}
&&\log(\Psi(t,u,v,w_1, w_2, \ldots, w_n)) + \int_0^t \sum_{i=1}^nh_2(y,w_i)dy\\
&= &
\int_0^t \sum_{i=1}^nh_2(y,w_i) \int_{{\mathbb R}^n}
h_4(t;u,v;  z_1,z_2, \ldots, z_n)\\ &&\Psi(t-y,u,v,z_1,z_2, \ldots, z_n)dz_1dz_2\ldots dz_ndy.
\end{eqnarray*}

In general, we have $f(x)$ in the following form $f(x)
=\sum_{\ell=0}^n f_{\ell}x^{\ell}$. Notice that if we replace
\begin{eqnarray*}
\exp\left[\sum_{k=1}^nh_1(t;u,v;z_k)\right],
\end{eqnarray*}
in (\ref{eqn:expan11}) by
\begin{eqnarray*}
\exp\left[\sum_{k=1}^\ell h_1(t;u,v;z_k)\right] \prod_{k=\ell+1}^n \nu(z_k).
\end{eqnarray*}
then the above logic applies to $x^\ell$, $\ell < n$, since $\nu$ has no
mass in $[0,M]$. Overall, we can apply the same procedure for general
function form $f(x)$, in which $h_4(t;u,v;
z_1,z_2, \ldots, z_n)$ takes the following form,
\begin{eqnarray*}
&&h_4(t;u,v;  z_1,z_2, \ldots, z_n) \\ &=& \sum_{\ell=0}^n
f_\ell \exp\left[\sum_{k=1}^\ell h_1(t;u,v;z_k)\right]
\prod_{k=\ell+1}^n \nu(z_k)\prod_{k=1}^n g(z_k).
\end{eqnarray*}
Therefore, we can rewrite the equation as,
\begin{eqnarray*}
&&\log(\psi(t,u,v,w))+  \int_0^t \sum_{i=1}^nh_2(y,w_i)dy \\
  &=&\int_0^t h_2(y,w)\left[\int_{{\mathbb R}^n}
\sum_{\ell=0}^n f_\ell \exp\left(\sum_{k=1}^\ell h_1(t;u,v;z_k)\right)
  \right. \\ && \left.\prod_{k=\ell+1}^n \nu(z_k)\prod_{k=1}^n
  g(z_k)\prod_{i=1}^n\psi(t-y,u,v,z_i)dz_1dz_2\ldots dz_n\right]dy ,
\end{eqnarray*}
with the convention that $\prod_{k=n+1}^n \nu(z_k) =1$. The lemma can
then be concluded in conjunction with Lemma \ref{lem:transform}.
$\Box$

\begin{thm}
\label{thm:solution}
{\rm The solution to the integral equation (\ref{eqn1_branch_2}) bears the following form,
\begin{eqnarray*}
\phi(t;u,v)& =& T \left(\exp[ h_4(t,u,v, z,z, \ldots,
z) -  n\int_0^t h_2(y,z)\right.\\ & &  \left. -
\int_0^t R(t-y, z,z,\ldots,z)   nh_2(y,z)]^\frac{1}{n}\right)
\end{eqnarray*}}
where $R(t, w_1, w_2, \ldots, w_n)$ is the inverse Laplace transform, with respect to $p$, of
\begin{eqnarray*}
\frac{\sum_{i=1}^n {\hat h}(p, w_i)}{1+\sum_{i=1}^n {\hat h}(p, w_i)},
\end{eqnarray*}
and ${\hat h}(p, w_i)$ denotes the Laplace transform of $h_2(t,w_i)$,
\begin{eqnarray*}
{\hat h}_2(p, w_i)= \int_0^\infty h_2(t,w_i) e^{-pt}dt.
\end{eqnarray*}
\end{thm}
{\bf Proof }
The equation (\ref{eqn:inteqn_n}) is a combined first and second kind integral
equation. The techniques for solving this type of equations are
of independent interest. For the one considered in the paper, the
solution to its component of integral equation of the first kind has a
special functional structure which enables us to solve the two components
separately.

According to \cite{polyaninmanzhirov}, the solution to the integral equation,
\begin{eqnarray*}
&& \log(\Psi(t,u,v,w_1, w_2, \ldots, w_n)) \\
& = & \int_{{\mathbb R}^n}
h_4(t,u,v, z_1,z_2, \ldots, z_n)\\&& \Psi(t,u, v, z_1, z_2, \ldots,
z_n))dz+\eta(t,u,v,w_1, w_2, \ldots, w_n)
\end{eqnarray*}
is in the following form,
\begin{eqnarray}
\label{eqn:additive}
&& \log(\Psi(t, u, v, w_1, w_2, \ldots, w_n)) \\& =& h_4(t,u,v,
w_1,w_2, \ldots, w_n) +\eta(t, u, v, w_1, w_2, \ldots, w_n). \nonumber
\end{eqnarray}
So it can be verified that if we can find a function $$\eta(t,u,v,w_1, w_2, \ldots,
w_n)$$ that satisfies,
\begin{eqnarray}
\label{eqn:linear}
&&\eta(t,u,v,w_1, w_2, \ldots, w_n) +  \sum_{i=1}^n\int_0^t h_2(y, w_i)\\
&=&\int_0^t   \sum_{i=1}^n h_2(y, w_i)\eta(t-y,u,v,w_1, w_2, \ldots, w_n)dy, \nonumber
\end{eqnarray}
and plug it into (\ref{eqn:additive}), $\Psi(t, u, v, w_1, w_2,
\ldots, w_n)$ will be the solution to (\ref{eqn:inteqn_n}). Meanwhile,
the integral equation (\ref{eqn:linear}) is a typical integral
equation of the second kind, therefore, can be solved routinely
by Laplace transform method. More specifically, we have, see e.g. \cite{BellmanCooke},
\begin{eqnarray*}
&&\eta(t,u,v,w_1, w_2, \ldots, w_n)= -  \sum_{i=1}^n\int_0^t h_2(y, w_i)\\&& -
\int_0^t R(t-y, w_1, w_2, \ldots, w_n)  \sum_{i=1}^n h_2(y, w_i) dy
\end{eqnarray*}
In conjunction with lemmas \ref{lem:transform} and \ref{lem:log_trans}, the result follows.
$\Box$

\section{Transform Expression of the Performance Measures and Inversion}
\label{sec:performance}

\vskip 0.2cm

In this section, we will first produce the final expression for the
Laplace transform of our stationary performance metric, based upon
results in the previous sections. Then, we will apply sophisticated
Laplace inversion methods to several well-known scheduling policies that
are either a special case of Grishechkin processor-sharing policy or
can be treated as its limiting case. Especially, following the approach in \cite{grishechkin}, we demonstrate that some popular scheduling policies such as foreground-backgound, SRPT and time-function scheduling can be expressed as the limit of a sequence Grishechkin processor-sharing policies, then, by carefully selecting the parameters, we are able to obtain the expressions for the stationary performance metrics and their Laplace transform in fairly tractable forms.

Given the solution of the integral equation, we can calculate the
Laplace transform of the performance metric $Q$.

\begin{thm}
{\rm The Laplace transform of the stationary system size is given by
\begin{eqnarray}
\label{eqn:laplace_final}
\ex\exp(-u Q) & = & (1-\rho) +
\int_0^\infty f'(\kappa_1(t,u,0)) \kappa_2(t,u)dt,
\end{eqnarray}
where
\begin{eqnarray*}
\kappa_1(t, u,v)&=& \int_{\mathbb R} \exp[h_1(t;u,v;z)]\times
\\&& \left(\exp[ h_4(t,u,v, z,z, \ldots,
z) -  n\int_0^t h_2(y,z)\right. \\ && \left. -
\int_0^t R(t-y, z,z, \ldots,
z)   n h_2(y,z)]^\frac{1}{n}\right)
g(z) dz
\end{eqnarray*}
\begin{eqnarray*}
\kappa_2(t,u) &= &  \int_{\mathbb R} \exp[h_1(t;u,0;z)] \times
\\&&  \left(\exp[ h_4(t,u,0, z,z, \ldots,
z) -  n\int_0^t h_2(y,z)  \right. \\ && \left.
- \int_0^t R(t-y, z,z, \ldots,
z)   n h_2(y,z)]^\frac{1}{n}\right)
\times \\&&  \frac{\partial}{\partial v} |_{v=0} h_1(t,u,v,z) g(z) dz
\\  &&+ \int_{\mathbb R} \exp[h_1(t;u,0;z)]\times \\&& \left(\exp[
  h_4(t,u,0, z,z, \ldots,z)  -  n\int_0^t h_2(y,z) \right. \\ &&
  \left. - \int_0^t R(t-y, z,z, \ldots,
z)   n h_2(y,z)]^\frac{1}{n}\right)\times \\&&
\frac{\partial}{\partial v} |_{v=0} \kappa_3(t,u,v,z) g(z)dz,
\end{eqnarray*}
and
\begin{eqnarray*}
\kappa_3(t,u,v,z)&=& [ h_4(t,u,0, z,z, \ldots,z)  -   n\int_0^t
  h_2(y,z) \\ &&-\int_0^t R(t-y, z,z, \ldots,
z)   n h_2(y,z)]^\frac{1}{n}.
\end{eqnarray*}}
\end{thm}

Apparently, this Laplace transform is in a very complicated functional
and integral form, it is unrealistic to seek its inversion in closed-form. Especially, the function $R(t,w_1,w_2,\ldots, w_n)$
is essentially an inverse Laplace transform, which makes our task essentially inverting a two-dimensional Laplace transform. This encourages us to use
numerical procedures for inverting Laplace transform to approximate
these functions. Extensive studies on a unified approach for inverting a Laplace transform are carried out in \cite{abatewhitt}. Among the methods discussed in \cite{abatewhitt}, we select the Talbot method for its concise expression and high accuracy. In the following we summarize the main idea of this method, for details, see, e.g. \cite{AbateValko} and \cite{abatewhitt}.

For any function $f$, the Laplace transform is defined by,
\begin{eqnarray*}
{\hat f}(s) = \int_0^\infty e^{-st} f(t) dt.
\end{eqnarray*}
Its inversion is given by the following Bromwich inversion integral,
\begin{eqnarray}
\label{eqn:Bromwich}
f(t) = \frac{1}{2\pi \sqrt{-1}} \int_C f(s) e^{st} ds, t> 0,
\end{eqnarray}
where the contour $C$ goes from $c-\infty \sqrt{-1}$ to $c+\infty
\sqrt{-1}$ for $c>0$. A unified Laplace inversion approach is to use
rational functions to approximate the exponential function in the
integrand. More specifically, use,
\begin{eqnarray*}
e^z \approx \sum_{k=0}^n \frac{w_k}{\al_k -z},
\end{eqnarray*}
for some carefully selected complex numbers $w_k$ and $\al_k$.
Then the Residue theorem will give a finite summation in
terms of the evaluation of Laplace transform for the Bromwich
integral (\ref{eqn:Bromwich}), more specifically,
\begin{eqnarray*}
f(t)= \frac{1}{t}\sum_{k=0}^n w_k {\hat f}\left(\frac{\al_k}{t}\right),t>0
\end{eqnarray*}
Different algorithms for Laplace inversion, such as the
Gaver-Stehfest algorithm, Euler algorithm and Talbot algorithm,
differ at the selection of the rational functions, i.e. $w_k$ and
$\al_k$. Here, we will use the Talbot algorithm.
Detailed analysis of the algorithm can be found in
\cite{abatewhitt, AbateValko}. For any large
integer $I>0$, the Talbot method uses the following expression as an
inversion of the Laplace transform ${\hat f}$.
\begin{eqnarray*}
f(t, I) = \frac{2}{5t} \sum_{k=0}^{I-1} {\rm Re} \left(\gamma_k {\hat
  f}\left(\frac{\delta_k}{t}\right)\right).
\end{eqnarray*}
with
\begin{eqnarray*}
\delta_0& = & \frac{2m}{5}, \\ \delta_k&=&
\frac{2k\pi}{5}\left[\cot\left(\frac{k\pi}{M}\right)+\sqrt{-1}\right], k>0\\
\gamma_0&=&\frac{1}{2}e^{\delta_0},\\ \gamma_k&= & \left\{1+
\sqrt{-1}\left(\frac{k\pi}{I}\right)\left[+\cot\left(\frac{k\pi}{I}\right)^2\right]\right.
\\ && \left.-\sqrt{-1}\cot
\left(\frac{k\pi}{I}\right)\right\}e^{\delta_k}, k>0.
\end{eqnarray*}

\subsection{Accuracy of the Laplace inversion}

\vskip 0.2cm

Here, we have a brief discussion on the accuracy of the Laplace
inversion so that we will have a full picture on the approach we are
taking.
\begin{defn}
For a large integer $M'>0$, and $\al >0$, we say that the inversion $f_q$ of
the Laplace transform ${\hat f}$ produces $\al M$ significant digits, if
\begin{eqnarray*}
\left| \frac{f(t)-f_q}{f(t)}\right| \approx 10^{-\al M'}
\end{eqnarray*}
\end{defn}
Then, it is known from \cite{AbateValko}, $f(t, I)$  the output of the Talbot inversion produces
$0.6I$ significant digits while requiring $I$ evaluations of the Laplace
transform. In the case of two-dimensional inverting Laplace transform,
such as the problem we are studying, it is demonstrated in
\cite{abatewhitt} that we can apply the Talbot algorithm to both, and
the overall algorithm still produces $0.6I$ significant digits while
requiring $I^2$ evaluations of the Laplace transform.

\subsection{Egalitarian Processor-sharing Queues}
\label{sec:PS}

\vskip 0.2cm

Let us consider the simplest egalitarian processor-sharing queue $M^X/G/1/PS$. In
this system, $C(t) =A(t) = {\bf 1} (t\in[0, \ell)]$. So $h_2(y, z) =
{\bf 1}\{ z \in [0, y]\}$, and
\begin{eqnarray*}
{\hat h}_2(p, w ) = \int_0^\infty  {\bf 1}\{ w \in [0, t]\}
e^{-pt}dt =\int_w^\infty e^{-pt} dt = e^{-pw}.
\end{eqnarray*}
Applying the Talbot method for integer $I_1>0$
\begin{eqnarray*}
&& R_1(t, w_1, w_2, \ldots, w_n, I_1)\\ &=& \frac{2}{5t}
  \sum_{k=0}^{I_1-1} {\rm Re} \left(\gamma_k
  \sum_{i=1}^n\frac{\exp(-w_i\delta_k/t)}{1+\sum_{i=1}^n\exp(-w_i\delta_k/t)}\right).
\end{eqnarray*}
Then the density function is given as, for integer $I_2>0$,
\begin{eqnarray}
\label{eqn:inverse_general}
&& \th_{Q}(s) = \frac{2}{5s} \sum_{k=0}^{I_2-1} {\rm Re} L_k,
\end{eqnarray}
where
\begin{eqnarray*}
L_k=
\gamma_k
\left[(1-\rho) +
\int_0^\infty f'(\kappa_1(t,\delta_k/t,0)) \kappa_2(t,\delta_k/t)dt\right] ,
\end{eqnarray*}
where $$R_1(t,w_1, w_2, \ldots, w_n,I_1)$$ will replace $$R(t,w_1, w_2,
\ldots, w_n)$$ in the definition of $\kappa_i, i=1,2,3$. If we let both
$I_1$ and $I_2$ be the order of some large integer $I>0$, then the
above expression produces  $0.6I$ significant digits.

In \cite{grishechkin}, for the case of egalitarian processor-sharing,
an integral equation is developed for another performance metric, the
sojourn time. More specifically, let $W(T)$ be the sojourn time for
a tagged job with processing time $\ell$, then, Theorem 6.2 in
\cite{grishechkin} gives,

\begin{thm} If $\rho< 1$ then $T\rightarrow \infty$, then $W(T)
  \rightarrow W$, and the Laplace transform of the random variable $W$
  is given by
\begin{eqnarray}
\label{eqn1_branch_4} &&\ex[\exp(-u W)] =  K(u, \ell)
  \exp(-u\ell) \exp\left[ \G \int_0^\ell (f(S(y,u))-1)dy
 \right],
\end{eqnarray}
where $f(\cdot)\in C^1$ is the Z-transform for the arrival batch
size. Here, $S(t,u)$ satisfies the following equation,
\begin{eqnarray}
\label{eqn:S_inteqn}
&&S(t,u)= \ex\exp\left[-u\min(t,\ell)-\G\int_0^\ell (1-f(S(t-y,u)))
  dy\right],
\end{eqnarray}
and
\begin{eqnarray*}
K(u,b) =\G(1-\rho)\left[\G^{-1} + \frac{\partial}{\partial v }\Big|_{v=0}
  \int_0^\infty f(\phi_b(t;u,v))dt\right]
\end{eqnarray*}
with $\phi_b$ satisfies,
\begin{eqnarray}
\label{eqn:pb_inteqn}
&&\phi_b(t;v,u) = \ex \exp\left[ - u
  \int_t^{t+b} {\bf 1}\{ y\in [0,\ell]\} dy \right.
\\&& \left. -v{\bf 1}\{ t\in [0,\ell]\} - \G\int_0^\ell
  1-f(\phi_b(t-y;v, u)) dy\right] \nonumber
\end{eqnarray}
\end{thm}
It should be clear that both $S(t,u)$ and $\phi(t,u,v)$ are in the
same form of integral equation as  (\ref{eqn1_branch_2}), hence,
similar techniques can be applied to them. More precisely, Theorem
\ref{thm:solution} can be applied to these two integral equations with
functions $h_1(\cdot)$ and $h_2(\cdot)$ defined as the following.
For  $S(t,u)$
\begin{eqnarray}
\label{eqn:H1_sojourn_S}
&&h_1(t;u,v;z_1,z_2,z_3,z_4)=[-u\min(t,\ell)(t)]{\bf 1}\{t=z_4\},
\end{eqnarray}
and
\begin{eqnarray}
\label{eqn:H2_sourjour_S}
h_2(t,z_1,z_2,z_3,z_4)={\bf 1}\{t=z_4\}.
\end{eqnarray}
and for $\phi_b$,
\begin{eqnarray}
\label{eqn:H1_sojourn_Pb}
&&h_1(t;u,v;z_1,z_2,z_3,z_4)\\&=&[ - u
  \int_t^{t+b} {\bf 1}\{ y\in [0,\ell]\} dy  -v{\bf 1}\{ t\in [0,\ell]\}]{\bf 1}\{t=z_4\}, \nonumber
\end{eqnarray}
and
\begin{eqnarray}
\label{eqn:H2_sourjour_Pb}
h_2(t,z_1,z_2,z_3,z_4)=C(z_1,z_2,z_3)(t){\bf 1}\{t=z_4\}.
\end{eqnarray}

After obtaining the solution of the integral equations, the numerical
procedure described above can be applied to get approximation with
the same performance guarantee.

\subsection{Discriminatory Processor-sharing Queues with random class assignment}
\label{sec:DPS}

\vskip 0.2cm

Discriminatory processor-sharing queue was first studied by Kleinrock
\cite{Kleinrock67}. Jobs are grouped in $C$ classes, indexed by
$c=1,2,\ldots, C$, each class carries a fixed weight $\nu_c$. In a
system with $N_c$ jobs for each class $c$ jobs, the amount of service
each job $x$ receives is determined by
\begin{eqnarray*}
\frac{ \nu_{c(x)}}{\sum_{c=1}^C N_c \nu_c} .
\end{eqnarray*}
where $c(x)$ denotes the job class $x$ belongs to. Of course, it is
easy to see that when all the $\nu_c$ are equal, the policy is just
the ordinary processor-sharing. For an updated survey on the analysis of
discriminatory processor-sharing policy, see,
e.g. \cite{AltmanAvrachenkovAyesta}.

Here, we consider a queue under a policy that is an variation of the
discriminatory processor-sharing policy. We allow the class
characterization being randomly determined at the time of
arrival. More specifically, a job will belong to class $c,
c=1,2,\ldots, C$ with probability
\begin{eqnarray*}
\frac{ \nu_c}{\sum_{c=1}^C\nu_c}.
\end{eqnarray*}
Therefore,
\begin{eqnarray*}
A(t) = \nu_c  {\bf 1}\{ w \in [0, t]\}
\end{eqnarray*}
and $C(t) =A(t) = {\bf 1} (t\in[0, \ell)]$. So $h_2(y, z) =
{\bf 1}\{ z \in [0, y]\}$, and
\begin{eqnarray*}
{\hat h}_2(p, w ) = \int_0^\infty \nu_c  {\bf 1}\{ w \in [0, t]\}
e^{-pt}dt =\int_w^\infty \nu_c e^{-pt} dt =\nu_c  e^{-pw}.
\end{eqnarray*}
Following our approach for the egalitarian processor-sharing, apply 
the Talbot method for integer $I_1>0$, we obtain, 
\begin{eqnarray*}
&& R_2(t, w_1, w_2, \ldots, w_n, I_1)\\ &=& \frac{2}{5t}
  \sum_{k=0}^{I_1-1} {\rm Re} \left(\gamma_k
  \sum_{i=1}^n\frac{\mu_c\exp(-w_i\delta_k/t)}{1+\sum_{i=1}^n\mu_c\exp(-w_i\delta_k/t)}\right).
\end{eqnarray*}
Plug the above $$R_2(t,w_1, w_2, \ldots, w_n,I_1)$$ as function $$R(t,w_1, w_2,
\ldots, w_n)$$ in the definition of $\kappa_i, i=1,2,3$ in equation (\ref{eqn:inverse_general}), for some $I_1$,  we obtain a guaranteed approximation for the density function of the performance depending on the selection of $I_1$ and $I_2$.

\subsection{Shortest Residual Processing Time First}
\label{sec:SRPT}

\vskip 0.2cm

In \cite{grishechkin}, it is shown that another popular queue
scheduling policy, the shortest residual processing time (SRPT) rule,
can be treated as the limit of a sequence of  Grishechkin
processor-sharing policies with
carefully chosen parameters. In particular, for any positive integer
$N=1,2,\ldots$, let $c_N$ be a sequence of functions satisfying
$c_N(T_1) /c_N(T_2)\rightarrow \infty$, as $N\rightarrow \infty$ for
any fixed $T_2>T_1>0$. Then define,
\begin{eqnarray*}
 A_N(T) = c_n(T) {\bf 1}\{T \le \ell\}.
\end{eqnarray*}
Denote $Q^\phi_N$ as the state process of the
$N$-th system, and  $Q^\phi$ as the state process of a system that
follows the SRPT rule, then Lemma 7.1 in \cite{grishechkin} indicates
that $Q^\phi_N$ converge to $Q^\phi$ in distribution, as $N\rightarrow
\infty$.

To facilitate the calculation, we select $c_N(y) = (\ell-y)^N$, it is easy
to see that this sequence satisfies the assumptions. Now,
\begin{eqnarray*}
 A_N(T) =  (\ell - T)^{N+1}{\bf 1}\{T \le \ell\}.
\end{eqnarray*}
\begin{eqnarray*}
R_N(u) = -\frac{1}{N} ( \ell^{-N} - (\ell-u)^{-N}
\end{eqnarray*}
\begin{eqnarray*}
C^N(u) = -\frac{1}{N}(\ell^N-Ny)^{-\frac{1}{N}-1}.
\end{eqnarray*}
Hence,
\begin{eqnarray*}
h_2(p,w) =  -\frac{1}{N} p^{-\frac{1}{N+1}}\left[ \gamma\left(\frac{1}{N+1},
\frac{\ell^N-w}{p}\right)-\gamma\left(\frac{1}{N+1},
\frac{\ell^N}{p}\right)\right].
\end{eqnarray*}
where $\gamma(s,x)$ denotes the lower incomplete Gamma function, which
is defined as $\gamma(s,x)= \int_0^xt^{s-1}e^{-t}dt$.
For any $M_1>0$, we have,
\begin{eqnarray*}
&& R_3^N(t, w_1, w_2, \ldots, w_n, M_1)= \frac{2}{5t} \sum_{k=0}^{M_1-1} {\rm Re} (Q_k),
\end{eqnarray*}
where $Q_k$ is defined in (\ref{eqn:Q}).
\begin{figure*}
\begin{eqnarray}
\label{eqn:Q}
Q_k=-N\frac{1-\frac{1}{N}\sum_{i=1}^n\gamma_k
  \frac{\delta_k}{t}^{-\frac{1}{N+1}}
  \left(\gamma\left(\frac{1}{N+1},
  \frac{t(\ell^N-w_i)}{\delta_k}\right)-\gamma\left(\frac{1}{N+1},
  \frac{t(\ell^N)}{\delta_k}\right)\right)}{\sum_{i=1}^n\gamma_k
  \frac{\delta_k}{t}^{-\frac{1}{N+1}}
  \left(\gamma\left(\frac{1}{N+1},
  \frac{t(\ell^N-w_i)}{\delta_k}\right)-\gamma\left(\frac{1}{N+1},
  \frac{t(\ell^N)}{\delta_k}\right)\right)}.
\end{eqnarray}
\end{figure*}

Again, plug the above $$R_3^N(t,w_1, w_2, \ldots, w_n,I_1)$$ as function $$R(t,w_1, w_2, \ldots, w_n)$$ in the definition of $\kappa_i, i=1,2,3$ in equation (\ref{eqn:inverse_general}) for properly chosen $I_2$, the density function can be computed to desired accuracy.


\subsection{The Foreground-Background Queue}
\label{sec:FB}

\vskip 0.2cm

Foreground-background policy is another policy that can be analyzed
using the techniques discussed in this paper, because it is a policy
that allocates the service capacity according to the service attained
for each individual job. To be more specific, foreground-background
policy gives priority to that job that has the least amount of service
attained, and if there is more than one such jobs, then the server is
shared equally among these jobs. It is known that under heavy-tailed
distribution, there is a strong positive correlation between large
attained service and large remaining service, therefore, under this
circumstance, foreground-background policy is considered as a
surrogate for SRPT when the total amount of service requirement for each job
can not be determined by the scheduler.

As pointed out in \cite{grishechkin}, similar logic to the above in
the case of the shortest remaining service rule applies here.
So, for any positive integer
$N=1,2,\ldots$, let $c_N$ be a sequence of functions satisfying
$c_N(T_1) /c_N(T_2)\rightarrow \infty$, as $N\rightarrow \infty$ for
any fixed $T_2>T_1>0$.
 Then define,
\begin{eqnarray*}
A_N(T) = c_N(T) {\bf 1}\{T \le \ell\}.
\end{eqnarray*}
Then, $Q^\phi_N$, defined  as the state process of the $N$-th system,
and  $Q^\phi$,  the state process of a system that
follows the FB rule, satisfy that $Q^\phi_N$ converge to $Q^\phi$ in
distribution, as $N\rightarrow \infty$.

To facilitate the calculation, we select $c_N(y) = y^{-N}$, it is easy
to see that this sequence satisfies the assumptions. From the
definition, we know that,
\begin{eqnarray*}
R_N(u) = \frac{1}{N+1}u^{N+1}, R^{-1}_N(U) = (N+1)^{\frac{1}{N+1}}u^{\frac{1}{N+1}}.
\end{eqnarray*}
and
\begin{eqnarray*}
C^N(u) = (N+1)^{\frac{-N}{N+1}}u^{\frac{-N}{N+1}}.
\end{eqnarray*}
Therefore,
\begin{eqnarray*}
{\hat h}_2(p,w)& = &(N+1)^{\frac{-N}{N+1}} \int_0^w
y^{\frac{-N}{N+1}}e^{-yp}dy\\
&=&(N+1)^{\frac{-N}{N+1}}p^{-\frac{1}{N+1}} \gamma\left(\frac{1}{N+1},
\frac{w}{p}\right).
\end{eqnarray*}
where $\gamma(s,x)$ denotes the lower incomplete Gamma function, which
is defined as $\gamma(s,x)= \int_0^xt^{s-1}e^{-t}dt$.
For any $I_1>0$, we have,
\begin{eqnarray*}
&& R_4^N(t, w_1, w_2, \ldots, w_n, M_1)= \frac{2}{5t} \sum_{k=0}^{I_1-1} {\rm Re} (Q_k),
\end{eqnarray*}
where
\begin{eqnarray*}
Q_k=\frac{1+\sum_{i=1}^n\gamma_k
  (N+1)^{\frac{-N}{N+1}}\frac{\delta_k}{t}^{-\frac{1}{N+1}}
  \gamma\left(\frac{1}{N+1},
  \frac{w_it}{\delta_k}\right)}{\sum_{i=1}^n\gamma_k
  (N+1)^{\frac{-N}{N+1}}\frac{\delta_k}{t}^{-\frac{1}{N+1}}
  \gamma\left(\frac{1}{N+1}, \frac{w_it}{\delta_k}\right)}.
\end{eqnarray*}
Plug the above $$R_4^N(t,w_1, w_2, \ldots, w_n,I_1)$$ as function $$R(t,w_1, w_2,
\ldots, w_n)$$ in the definition of $\kappa_i, i=1,2,3$ in equation (\ref{eqn:inverse_general}), we obtain a guaranteed approximation for the density function of the performance metric.

\subsection{Time Function Scheduling}
\label{sec:TFS}

\vskip 0.2cm

Time function scheduling is another scheduling policy that was first
studied by Kleinrock, see, e.g. \cite{Kleinrock}. Further studies can
be found in \cite{FongSquillante} and \cite{LuSquillante}. Under this policy,
jobs are grouped in $C$ classes, and a weight $\nu_c>0$ is
assigned to each class $c=1,2,\ldots, C$. The server serves only one
job at a time. For any class $c$ job in the system, a time function is
calculated in terms of its cumulative waiting time and its class
weight $\nu_c$, in the literature, such as, \cite{Kleinrock},
\cite{FongSquillante} and \cite{LuSquillante}, this function is
basically taking the form of a linear function of the waiting time with
$\nu_c$ as the slope. The scheduling policy is to assign the job that
has the highest value of its time function to be served by the
server. Again as in Sec. \ref{sec:DPS}, we assume that a job is
randomly assigned to a class, and the probability distribution is denoted by
$\mu_c$.

Following similar derivations as those for the SRPT rule in
\cite{grishechkin}, we can show that the above defined time function
scheduling policy can also be treated as a limit of Grishechkin
processor-sharing policies. For any positive integer
$N=1,2,\ldots$, let $c_N$ be a sequence of functions satisfying
$c_N(T_1) /c_N(T_2)\rightarrow \infty$, as $N\rightarrow \infty$ for
any fixed $T_2>T_1>0$. Then define,
\begin{eqnarray*}
 A_N(T) = c_N(T) \mu_c {\bf 1}\{T \le \ell\}.
\end{eqnarray*}
Now, denote $Q^\phi_N$ as the state process of the
$N$-th system following a Grishechkin processor-sharing policy with
$A$ defined as above, and  $Q^\phi$ as the state process of a system
that follows the time function rule. We can demonstrate that
$Q^\phi_N$ converge to $Q^\phi$ in distribution, as $N\rightarrow
\infty$.

Now, let us select $c_N(y) = y^{-N}$. Following the same calculation as in foreground-background queue, we have,  Hence,
\begin{eqnarray*}
h_2(p,w)=\nu_c(N+1)^{\frac{-N}{N+1}}p^{-\frac{1}{N+1}} \gamma\left(\frac{1}{N+1},
\frac{w}{p}\right).
\end{eqnarray*}
Similarly, for any $I_1>0$, we have,
\begin{eqnarray*}
&& R_5^N(t, w_1, w_2, \ldots, w_n, M_1)= \frac{2}{5t} \sum_{k=0}^{I_1-1} {\rm Re} (Q_k),
\end{eqnarray*}
where
\begin{eqnarray*}
Q_k=\frac{1+\sum_{i=1}^n\mu_c\gamma_k
  (N+1)^{\frac{-N}{N+1}}\frac{\delta_k}{t}^{-\frac{1}{N+1}}
  \gamma\left(\frac{1}{N+1},
  \frac{w_it}{\delta_k}\right)}{\sum_{i=1}^n\mu_c\gamma_k
  (N+1)^{\frac{-N}{N+1}}\frac{\delta_k}{t}^{-\frac{1}{N+1}}
  \gamma\left(\frac{1}{N+1}, \frac{w_it}{\delta_k}\right)}.
\end{eqnarray*}
The rest of the calculation can follow the same approach as for the foreground-background policy.

\section{Conclusions}
\label{sec:conclusions}

\vskip 0.2cm

In this paper, we analyze the stationary performance of a queueing
system under a general class of processor sharing scheduling
policies. Our main contribution is obtaining a solution to a complicated integral equation that plays a critical role in queueing analysis. The methods we derived in solving the integral equation appear to be of independent interest to many other problems in mathematics and engineering. Meanwhile, we adopted a sophisticated numerical Laplace inversion scheme, so that the relative error of the numerical inversion can be easily controlled. These results have important implications in the development of numerical computational package softwares for performance analysis and optimal control.

As we have demonstrated in the paper, because they allow service capacities to be dynamically determined by the attained service for each job,  this general class of scheduling policies, Grishechkin processor-sharing policies, are very powerful and flexible mathematical models. Their performance analysis leads to deeper understanding and analysis of many popular scheduling policies. It also provides a building block for potential development of adaptive scheduling policies that can be used for different performance requirements. A typical example is to eliminate the independence assumptions on stochastic processes $A_i$ for different $i$. More precisely, let $(A_1, A_2, \dots, A_n)$ be dependent stochastic processes indexed by an infinite-dimensional sequence of service attained for each job, then,
\begin{eqnarray*}
\frac{dV_i(t)}{dt} = \frac{ A_i(V_1, V_2,\ldots, ) }{\sum_{j=1}^{Q(t)} A_j(V_1, V_2,\ldots, )}
\end{eqnarray*}
We aim to use some more generalized branching process to characterize this type of scheduling policy. Once this can be achieved, policies such as SRPT, FB and time function scheduling can be directly incorporated instead of resorting to limit.

The integral equation studied in the paper has a very complicated form, however, it is not extremely eccentric. Quite often, applications in mathematics, engineering and economics also produce integral equations  of a similar form. In fact, some stochastic control and optimal stopping problems
appear to be very closely related to integral equations of this type,
see, e.g. a recent result on evaluation of finite horizon
Russian option in  \cite{Peskir}. Therefore, our other line of
research is to apply some of the techniques we developed here to
integral equations arising in other areas.


\begin{thebibliography}{99}
\baselineskip 12pt

\vskip 0.3cm


\bibitem{AbateValko}
{\sc Abate, J. \ and Valko, P. P. } (2004) Multi-Precision Laplace inversion,
{\it Int. J. Num. Meth. Engng}, 60(5):979-993.


\bibitem{abatewhitt}
{\sc Abate, J. \ and Whitt, W. } (2006) A Unified Framework for Numerically
Inverting Laplace Transforms. {\it INFORMS Journal on Computing}, 18(4): 408-421.

\bibitem{AltmanAvrachenkovAyesta}
{\sc Altman, E., Avrachenkov, K. \ and Ayesta, U. } (2006) A survey on
discriminatory processor sharing, {\it Queueing Systems}, 53(1-2): 53-63.





\bibitem{asmussenbook}
{\sc Asmussen, S. } (1987) {\it Applied Probability and Queues}, Wiley,
Chichester.

\bibitem{AABN05}
{\sc  Avrachenkov, K. Ayesta, U., Brown, P. \ and Nunez-Queijia, R., } (2005)
Discriminatory processor sharing revisited, {\it INFOCOM}.

\bibitem{baccelli}
{\sc Baccelli, F. \ and Towsley, D. } (1990) The customer response times in
the processor sharing queue are associated, {\it Queueing Systems}, 7(3-4):
269-282.


\bibitem{BellmanCooke}
{\sc Bellman, R. \ and Cooke, K. L. ,} (1963) {\it Differential-Difference
  Equations} Academic Press, New York, .



\bibitem{FongSquillante}
{\sc Fong, L. \ and  Squillante, M., } (1995) Time-Function Scheduling: a
general approach to controllable resource management, {\it
  Proc. Symp. Op. Sys. Prin.} 29(5):230-230.

\bibitem{FriedmanHenderson}
{\sc Friedman, E. J. \ and Henderson, S. G. } (2003) Fairness and efficieny
in web server protocols, in {\it Proceedings of ACM/SIGMETRIC'03}.


\bibitem{grishechkin91a}
{\sc  Grishechkin, S. } (1991) Crump-Mode-Jagers Branching Processes as A
Method for Investigation the $M/G/1$ System with Processor Sharing,
{\it Teor. Verojatn. i. Primen},  36(1):16-33.


\bibitem{grishechkin}
{\sc Grishechkin, S. } (1992) On a Relationship Between Processor-Sharing
Queues and Crump-Mode-Jagers Branching Processes, {\it Advances in Applied
Probability } 24(3):653-698.


\bibitem{gromoll2002}
{\sc Gromoll, C., Puha, A \ and Williams, R., } (2002) The fluid limit of a
heavily loaded processor sharing queue, {\it Annals of Applied
  Probability}, 12(3):797-859.

\bibitem{gromoll2004}
{\sc Gromoll, C., } (2004) Diffusion approximation of a processor sharing
queue in heavy traffic. {\it Annals of Applied Probability}, 14(2):
555-611.


\bibitem{Jagers}
{\sc Jagers, P. } (1975) {\it Branching Processes with Biological
  Applications. } Wiley, Chichester.

\bibitem{JagersNerman}
{\sc Jager, P. \ and Nerman, O. } (1984) The growth and composition of
branching populations, {\it Advances in Applied Probability}, 16(2):221-259.


\bibitem{jean-marie}
{\sc Jean-Marie, A \ and Robert, P} (1994) On the transient behavior of the
processor sharing queue, {\it Queueing Systems}, 17(1-2): 129-136.

\bibitem{Kleinrock67}
{\sc Kleinrock, L}, (1967) Time-Shared Systems: a theoretical treatment,
JACM, 14(2):242-261.



\bibitem{Kleinrock}
{\sc Kleinrock, L}, (1976) {\it Queueing Systems, Vol II: Computer
  Applications, } Wiley-Interscience.

\bibitem{LuSquillante}
{\sc Lu, Y. \ and Squillante, M. } (2005) Dynamic Scheduling to minimize utility
functions of sojourn time moments in queueing systems, {\it
  Performance Evaluation Review}, 33(2):42-44, .

\bibitem{NuyensWierman}
{\sc Nuyens, M. \ and Wierman, A. } (2008) The Foreground-Background queue: A survey
{\it Performance Evaluation}, 65(3-4):289-307.

\bibitem{Peskir}
{\sc Peskir, G. } (2005) The Russian option: finite horizon, {\it Finance and
  Stochastics}, 9(2):251-267, .



\bibitem{polyaninmanzhirov}
{\sc Polyanin, A. D. \ and Manzhinov, A.V. } (1988) {\it Handbook of Integral
  Equations}, CRC Press.

\bibitem{PorterStirling}
{\sc Porter, D. \ and Stirling, D.} (1990) {\it Integral Equations},
Cambridge Texts in Applied Mathematics, Cambridge University Press.

\bibitem{Wierman}
{\sc Wierman, A. }(2010) Fairness and scheduling in single server queues,
{\it Surveys in Operations Research and Management Science}.

\bibitem{yashkov87}
{\sc Yashkov, S.F., } (1987) Processor-sharing queues: some progress in
analysis, {\it Queueing Systems}, 2(1-2):1-17.




\end{thebibliography}
\end{document}